\documentclass[11pt]{amsart}
\usepackage{latexsym,color,amsmath,amsthm,amssymb,amscd,amsfonts}
\usepackage{amssymb, amsmath, latexsym, amsthm, amscd, graphicx, xcolor, hyperref, bm,enumitem}
\usepackage{nicefrac}
\def\R{\mathbb R}
\def\N{\mathbb N}

\def\Z{\mathbb Z}
\def\Q{\mathbb Q}

\def\={\equiv}
\def\<{\langle}
\def\>{\rangle}

\def\1{\mathbf{1}}

\def\opt1{\operatorname{T}^1}
\def\SL{\operatorname{SL}}

\newcommand{\df}{{\, \stackrel{\mathrm{def}}{=}\, }}

\numberwithin{equation}{section}

\newcommand{\sm}{\smallsetminus}

\newtheorem{theorem}{Theorem}[section]
\newtheorem{proposition}[theorem]{Proposition}

\newtheorem{lemma}[theorem]{Lemma}

\newtheorem{remark}[theorem]{Remark}

\newtheorem{conjecture}[theorem]{Conjecture}

\title{On Non-Uniformly Discrete Orbits}
\author[]{Sahar Bashan}
\thanks{}

\begin{document}
\begin{abstract}
We study the property of uniform discreteness within discrete orbits of non-uniform lattices in 
$SL_2(\R)$, acting on $\R^2$ by linear transformations. We provide quantitative consequences of previous results by using Diophantine properties. We give a partial result toward a conjecture of Lelièvre regarding the set of long cylinder holonomy vectors of the "golden L" translation surface: for any $\epsilon>0$, three points of this set can be found on a horizontal line within 
a distance of $\epsilon$ of each other.
\end{abstract}
\maketitle

 \section{Introduction}
\begin{figure}[ht]
    \centering
    \includegraphics[scale=0.8]{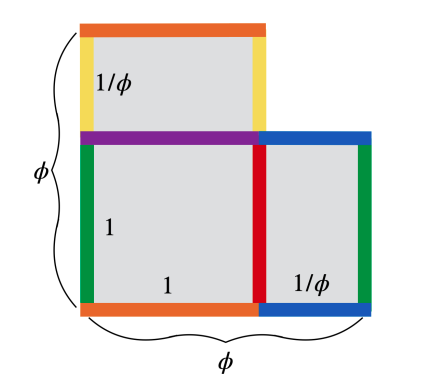}
    \caption{The golden L translation surface, with edge identifications.}
    \label{fig:GoldenL}
\end{figure}

A discrete set $Z \subset \R^2$ is called {\em uniformly discrete} if
$$
\inf\{ \| z_1 - z_2 \| : z_1, z_2 \in Z, \ z_1 \neq z_2 \} >0. $$
 Let $\Gamma=\bigr\langle\sigma_0,\sigma_1,\sigma_2,\sigma_3\bigr\rangle$ where
\[\sigma_0=
 \begin{pmatrix}
1 & \phi\\
 0 & 1
 \end{pmatrix} \qquad
 \sigma_1=
 \begin{pmatrix}
 \phi & \phi\\
 1 & \phi
 \end{pmatrix}\]

 \[\sigma_2=
 \begin{pmatrix}
 \phi & 1\\
 \phi & \phi
 \end{pmatrix}\qquad
 \sigma_3=
 \begin{pmatrix}
 1 & 0\\
 \phi & 1
 \end{pmatrix}\]
 and $\phi=\frac{1+\sqrt 5}{2}$ is the golden ratio.
 \begin{remark}
 The group $\Gamma$ is generated by the matrices $\sigma_0$ and $\sigma_3$ alone (see Lemma 2.7 and Remark 2.8 in \cite{Davis_Lelièvre}). However, for the purpose of the proofs and constructions in the following sections, we will make use of all four generators.
 \end{remark}
 In this paper we look at the discrete set $$S=\bigr\{\gamma\big(\begin{smallmatrix}
 1\\
 0
  \end{smallmatrix}\big):\gamma\in \Gamma\bigr\}.$$ It is known that $S$ is not uniformly discrete as it is the orbit of a point under the action of a non-arithmetic lattice in $SL_2(\R)$ (see \cite{Wu_thesis}). We show the following:
\begin{proposition}\label{prop}
 For every $0<\epsilon\leq1$, distinct points of $S$ in distance at most $\epsilon$ can be found in the ball $B(0,r)$ where $r=O(\frac{1}{\epsilon^2})$ .
\end{proposition}
The motivation for studying the set $S$ comes from the world of Veech surfaces.
$\Gamma$ is the Veech group of the golden L, which is a non-arithmetic  lattice surface (see figure \ref{fig:GoldenL}). The union of the orbits $$S\cup 
 \phi^{-1}S$$ is the set of holonomy vectors of the golden L (see \cite{Davis_Lelièvre}). For information on Veech surfaces see \cite{Smillie_Weiss} and \cite{Veech}.
 
In \cite{Wu}, Chenxi Wu proved that the set of holonomy vectors of a non-arithmetic Veech surface is not uniformly discrete. As mentioned, he also proved the generalization for orbits of any point in $ \R^2\setminus\{0\}$ under a non-arithmetic lattice in $SL_2(\R)$.  
Here we recall Wu's argument and add a quantitative consequence. As in \cite{Wu_thesis}, we use the fact that every non-uniform lattice in $\SL_2(\R)$ has 2 non-commuting unipotent elements (see \cite{MT}). We look at the group $\Gamma_\rho :=\langle \left( \begin{smallmatrix} 1 &
   1\\ 0 & 1 \end{smallmatrix} \right), \left( \begin{smallmatrix} 1 & 0 \\ \rho & 1 \end{smallmatrix} \right)\rangle$, for $\rho>0$.
   \begin{proposition}\label{rho}
$\Gamma_\rho\big(\begin{smallmatrix}
    1\\0
\end{smallmatrix}\big)$ is
uniformly discrete if and only if $\rho \in \mathbb{Q}$.
\end{proposition}
For the irrational case the argument will utilize the convergents of the continued fraction of the parameter $\rho$. The convergents will be used also for the asymptotics of the non-uniform discreteness:
   \begin{proposition}\label{con}
   Assume $\rho \notin\mathbb{Q}$ and let $\big(\frac{p_n}{q_n}\big)$ be the convergents of $\rho+1$. Then, in $ B(0, r) $, where $r=O\left(q_nq_{n+1}\right) $, there are distinct points of the orbit $\Gamma_\rho\big(\begin{smallmatrix}
    1\\0
\end{smallmatrix}\big)$ at a distance of $ O\left(\frac{1}{q_{n+1}}\right) $ from each other.
\end{proposition}
Thus, the Diophantine properties of $\rho$ appear in the asymptotic bounds. Recall that the continued fraction denominators satisfy the relation $q_{n+1}=a_{n+1}q_n+q_{n-1} \approx a_{n+1}q_n$, where $a_n$ are the continued fraction partial quotients of $\rho$ (they are the same for $\rho,\rho+1$, for $n>0$). If we construct $\rho$ by choosing the convergents, such that the partial quotients $a_n$ are sufficiently large at each step, then $\rho$ is very well approximated by rationals. As a result, according to Proposition \ref{con}, we can find distinct points in the orbit $\Gamma_\rho\big(\begin{smallmatrix}
    1\\0
\end{smallmatrix}\big)$ at a distance of $O(\epsilon_n)$ from each other, where $\epsilon_n=\frac{1}{q_{n+1}}$, in a ball of radius $O(\epsilon_n^{-1})$, centered in the origin.
 \\
Samuel Lelièvre had a stronger conjecture regarding the set $S$:
\begin{conjecture}\label{conjecture}
For every $\epsilon >0$ 
 and for every $n\in \N$ there are $x_1,...,x_n\in S$ on a horizontal line, such that $$\forall 1\leq i,j\leq n, \ \ \  \|x_i-x_j\|<\epsilon.$$ \end{conjecture} 
The main result of this paper is a proof of Conjecture \ref{conjecture} in case $n=3$. We prove the following theorem:  
  \begin{theorem}\label{3points}
 There exists a sequence of positive numbers $\epsilon_n\rightarrow0$ such that for every $n$, there are 3 points of $S$ in a horizontal line segment of length $\epsilon_n$. These points lie in a ball of radius $O(\epsilon_n^{-4})$, centered in the origin.
  \end{theorem}
The author wishes to thank Barak Weiss for his invaluable guidance and comments on earlier versions of this paper. The author also thanks Samuel Lelièvre and Julian Rüth for their helpful assistance and discussions, Shachar Heyman for his support with programming, and Claire Burrin for her valuable comments. Thanks are also due to the anonymous referee for useful comments. This paper is part of the author’s Master’s thesis at the School of Mathematics at Tel Aviv University, conducted under the supervision of Barak Weiss. The support of grant ISF-NSFC 3739/21 is gratefully acknowledged.
 \section{Uniform discreteness and continued fractions}
We denote by $\{\cdot\}$ the fractional part of a non-negative real number, i.e., $\{x\} \df x-\lfloor x \rfloor$.
For the proof of Proposition \ref{prop}, we utilize the technique from \cite{Wu_thesis}, adapting it to specific points. This approach provides a quantitative result regarding the growth rate of the norms of nearby orbit points. 
\begin{proof}[Proof of Proposition \ref{prop}]
    Let $(F_n)_{n\in \mathbb{N}}$ be the Fibonacci sequence, i.e., $$F_0=0, \hspace{1mm} F_1=1,$$ $$\forall n>1, \hspace{1mm}F_{n+1}=F_n+F_{n-1}.$$ 
    Notice that $$\phi^{-(n-1)}=(-1)^{n-1}\cdot F_n+(-1)^{n}\cdot F_{n-1}\phi,$$
   and for an odd positive integer $k$, $$\lfloor\phi^k\rfloor=F_{k-1}+F_{k+1}.$$ Additionally, $(\frac{F_{n+1}}{F_n})$ are the convergents to $\phi$. \\Let $n$ be an even positive integer. Define $\gamma_1,\gamma_2\in \Gamma$ by: \[
\gamma_1=\sigma_3^{F_{n-4}+F_{n-2}}\sigma_0^{F_n}\sigma_2,\]
\[\gamma_2=\sigma_3^{F_{n-4}+F_{n-2}}\sigma_0^{F_{n+1}}\sigma_1
.\]
We now compute as follows:
\[\big\|\gamma_1\big(\begin{smallmatrix}
 1\\
 0
  \end{smallmatrix}\big)-\gamma_2\big(\begin{smallmatrix}
 1\\
 0
\end{smallmatrix}\big)\big\|=\Big\|\sigma_3^{F_{n-4}+F_{n-2}}\Bigl(\big(\begin{smallmatrix}
 \phi+F_n\phi^2\\
 \phi
  \end{smallmatrix}\big)-\big(\begin{smallmatrix}
 \phi+F_{n+1}\phi\\
 1
  \end{smallmatrix}\big)\Bigr)\Big\|
  \]
  \[
  =\Big\|\sigma_3^{F_{n-4}+F_{n-2}}\big(\begin{smallmatrix}
 F_n-F_{n-1}\phi\\
 \phi-1
  \end{smallmatrix}\big)\Big\|
  =\Big\|\sigma_3^{F_{n-4}+F_{n-2}}\big(\begin{smallmatrix}
 -\phi^{-(n-1)}\\
 \phi^{-1}
  \end{smallmatrix}\big)\Big\|\]
  \[
  =\Big\|\Bigl(\begin{smallmatrix}
 -\phi^{-(n-1)}\\
 \phi^{-(n-2)}(\phi^{n-3}-(F_{n-4}+F_{n-2}))
\end{smallmatrix}\Bigr)\Big\|=\Big\|\Bigl(\begin{smallmatrix}
 -\phi^{-(n-1)}\\
 \phi^{-(n-2)}(\phi^{n-3}-\lfloor\phi^{n-3}\rfloor)
  \end{smallmatrix}\Bigr)\Big\|\]
  \\
  \[
  =\Big\|\Bigl(\begin{smallmatrix}
 -\phi^{-(n-1)}\\
 \phi^{-(n-2)}\{\phi^{n-3}\}
\end{smallmatrix}\Bigr)\Big\|\xrightarrow[n \to \infty]{}0.
\]
  \\
This computation shows that the distances between distinct points in S are not bounded below. Now, we bound the norm  $\big\|\gamma_2\big(\begin{smallmatrix}
 1\\
 0
\end{smallmatrix}\big)\big\|$:
$$\big\|\gamma_2\big(\begin{smallmatrix}
 1\\
 0
\end{smallmatrix}\big)\big\|=\Big\|\Bigl(\begin{smallmatrix}
 (F_{n+1}+1)\phi\\
 1+(F_{n-2}+F_{n-4})(F_{n+1}+1)\phi^2
\end{smallmatrix}\Bigr)\Big\|$$$$=\Big\|\Bigl(\begin{smallmatrix}
 (F_{n+1}+1)\phi\\
 1+(F_{n+1}+1)\lfloor\phi^{n-3}\rfloor\phi^2
\end{smallmatrix}\Bigr)\Big\|\leq(F_{n+1}+1)\phi+1+(F_{n+1}+1)\phi^{n-1}=O(\phi^{2n}).$$ \\
Here we used $F_n=\frac{1}{\sqrt5}(\phi^n-(-\phi)^{-n}).$\\
Let $0<\epsilon\leq1$, and $n$ such that $\phi^{-(n+1)}<\epsilon\leq \phi^{-n}$. We've shown that there are 2 distinct points of the orbit $\Gamma_\rho\big(\begin{smallmatrix}
    1\\0
\end{smallmatrix}\big)$ in distance $O(\phi^{-n})=O(\epsilon)$ from each other, in a ball of radius $O(\phi^{2n})=O(\epsilon^{-2})$.

\end{proof}
 Notice that we act with the unipotent matrix $\sigma_0$ in order to obtain two points with close first coordinate. In the definition of $\gamma_1,\gamma_2$, the powers of $\sigma_0$  are the convergents of the continued fraction of $\phi$. \\
We now generalize this idea to any finite union of discrete orbits of a non-uniform
lattice in $\SL_2(\R)$, acting on $\R^2$ by linear transformations. 
Recall that a matrix in $\SL_2(\R)$ is said to be \textit{unipotent} if its
eigenvalues are both equal to one. Note that if $\Gamma$ is a non-uniform lattice in $\SL_2(\R)$, that is, the quotient $\SL_2(\R)\big/\Gamma$ is not compact, then it
contains two non-commuting unipotent elements $u_1, u_2$. Also, it
follows from a theorem of Dani and Raghavan \cite{Dani_Raghavan}
that if $\Gamma$ is a
non-uniform lattice in $\SL_2(\R)$ and $z \in \R^2 \sm \{0\}$ is such
that the orbit $\Gamma z$ is discrete, then $\Gamma$ contains a
unipotent element fixing 
$z$. Up to a conjugation in
$\SL_2(\R)$, and up to possibly replacing $u_i$ with $u_i^{-1}$, we can assume that
$$
u_1 = \left( \begin{matrix} 1 &
   1\\ 0 & 1 \end{matrix} \right), \ \ \
\ u_2 = \left( \begin{matrix} 1 & 0 \\ \rho & 1 \end{matrix} \right), 
$$
for some $\rho = \mathrm{tr}(u_1u_2)-2> 0$. See \cite{KS} and the
references therein for 
information about the structure of the groups $\Gamma_{\rho} =
\langle u_1, u_2\rangle$.\\
Proposition \ref{con} demonstrates that $\rho$ also plays 
a role in quantifying the problem of non-uniform discreteness of $\Gamma_\rho
\mathbf{e}_1$.

\begin{proof}[Proof of Proposition \ref{con}] Consider the points:
$$\gamma_1 = u_1^L u_2^{p_n} u_2 \begin{pmatrix} 1 \\ 0 \end{pmatrix} \quad \text{and} \quad \gamma_2 = u_1^L u_2^{q_n} u_1 u_2 \begin{pmatrix} 1 \\ 0 \end{pmatrix},$$
 where $L=\bigl\lfloor(p_n-q_n(\rho+1))^{-1}\bigr\rfloor.$ Our proof proceeds in two steps: first, we bound the distance between these points; second, we estimate their norms. \\ WLOG we assume $p_n-q_n(\rho+1)>0$.\\
We compute the difference between 
$\gamma_1$ and $\gamma_2$:
\begin{align*}
\gamma_1 - \gamma_2 &= u_1^L \left( u_2^{p_n} u_2 \begin{pmatrix} 1 \\ 0 \end{pmatrix} - u_2^{q_n} u_1 u_2 \begin{pmatrix} 1 \\ 0 \end{pmatrix} \right) \\
&= u_1^L \left( \begin{pmatrix} 1 \\ (p_n+1)\rho \end{pmatrix} - \begin{pmatrix} \rho +1 \\ q_n\rho(\rho +1)+\rho \end{pmatrix} \right) \\
&= u_1^L \begin{pmatrix} -\rho \\ \rho(p_n - q_n(\rho +1)) \end{pmatrix} \\
&= \begin{pmatrix} -\rho + L \cdot \rho(p_n - q_n(\rho +1)) \\ \rho(p_n - q_n(\rho +1)) \end{pmatrix} \\
&= \rho \begin{pmatrix} -1 + L \cdot (p_n - q_n(\rho +1)) \\ p_n - q_n(\rho +1) \end{pmatrix}.
\end{align*}

Let $\Psi_n = p_n - q_n(\rho+1)$. Using the property $$L\cdot (p_n-q_n(\rho +1))=\Psi_n\lfloor\Psi_n^{-1}\rfloor=\Psi_n(\Psi_n^{-1}-\{\Psi_n^{-1}\})=1+O(\Psi_n),$$ we find that the difference simplifies to:
$$\rho \begin{pmatrix} O(\Psi_n) \\ \Psi_n \end{pmatrix}.$$  
Since $p_n$ and $q_n$ are convergents, we know that $|\Psi_n| = |p_n - q_n(\rho+1)| < \frac{1}{q_{n+1}}$ (see Theorem 9 in \cite{Khinchin}). Therefore, the distance between the points is bounded by $O\left(\frac{1}{q_{n+1}}\right)$.
 \\ 
 Next, we bound the norm of $\gamma_2 = u_1^L u_2^{q_n} u_1 u_2 \begin{pmatrix} 1 \\ 0 \end{pmatrix}$.
\begin{align*}
\gamma_2 &= u_1^L u_2^{q_n} \begin{pmatrix} 1+\rho \\ \rho \end{pmatrix} \\
&= u_1^L \begin{pmatrix} 1+\rho \\ \rho + \rho q_n(1+\rho) \end{pmatrix} \\
&= \begin{pmatrix} 1+\rho + L \rho (1+q_n+q_n\rho) \\ \rho + \rho q_n(1+\rho) \end{pmatrix}.
\end{align*}
The second coordinate is clearly $
O(q_n)$. For the first coordinate, we have:
\begin{align*}
&\left| 1+\rho + \rho L (1+q_n+q_n\rho) \right| \\
&\leq 1+\rho+\rho \left| (p_n - q_n(\rho+1))^{-1} \right| (1+q_n(\rho +1)) \\
&\leq c \left( q_n(q_{n+1}+q_n) \right).
\end{align*}
In the last inequality we used the property:
$$
\frac{1}{q_n(q_{n+1}+q_n)} < \left|\frac{p_n}{q_n}-(\rho+1)\right|
$$
(see Theorem 13 in \cite{Khinchin}).

Finally, combining these bounds, the norm of $\gamma_2$ is:
$$
\left\|\gamma_2\right\| \leq c'(q_n q_{n+1}+q_n^2+q_n) = O(q_n q_{n+1}).
$$
\end{proof}

\begin{proof}[Proof of Proposition \ref{rho}]
Proposition \ref{con} gives one direction. For the other one, if $\rho = \frac{p}{q} \in \mathbb{Q}$ then $\Gamma_\rho$ is contained in the
group
$$
\Lambda \df \left \langle 
 \left( \begin{matrix} 1 &
   1\\ 0 & 1 \end{matrix} \right), \ \left( \begin{matrix} 1 & 0 \\
   1/q & 1 \end{matrix} \right) 
\right
\rangle.
$$
All elements
of $\Lambda$ have entries which are rational numbers with
denominators which (in reduced form) divide $q$. This implies that
$$\Gamma_\rho\mathbf{e}_1 \subset \Lambda \mathbf{e}_1 \subset \frac{1}{q}
\, \Z^2,$$
and hence $\Gamma_\rho\big(\begin{smallmatrix}
    1\\0
\end{smallmatrix}\big)$ is uniformly discrete.
\end{proof}

\section{Weak Uniform discreteness}

\begin{figure}[ht]
    \centering
    \includegraphics[scale=0.4]{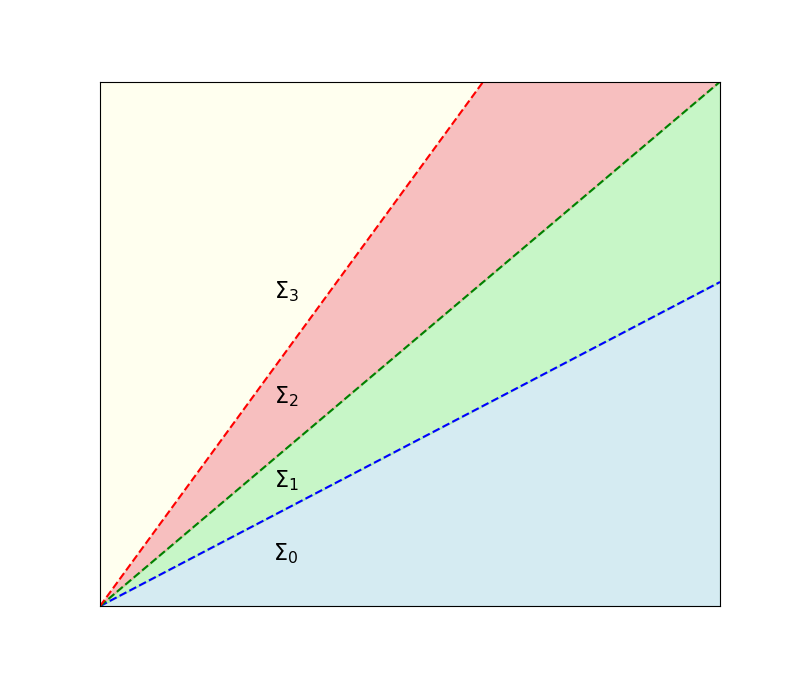}
    \caption{}
    \label{fig:partition}
\end{figure}

In \cite{Davis_Lelièvre}, Davis$\And$Lelièvre defined an iterative process for finding an element of S in the direction of a given vector v with slope in $\mathbb Q\bigr[\sqrt 5\bigr]$, in the first quadrant. The process is as follows:
let $\bigr\{\Sigma_0,\Sigma_1,\Sigma_2,\Sigma_3\bigr\}$ be a partition of the first quadrant (see Figure \ref{fig:partition}), such that
\[\Sigma_0=\bigr\{(x,y):x\geq0\hspace{1mm},\hspace{1mm}0\leq y<\phi^{-1}x\bigr\} \ \ \ \
\Sigma_1=\bigr\{(x,y):x\geq0\hspace{1mm},\hspace{1mm}\phi^{-1}x\leq y<x\bigr\}\]
\[\Sigma_2=\bigr\{(x,y):x\geq0\hspace{1mm},\hspace{1mm}x\leq y<\phi x\bigr\} \ \ \ \ \ \ \
\Sigma_3=\bigr\{(x,y):x\geq0\hspace{1mm},\hspace{1mm}\phi x\leq y\bigr\}.\]
\\
Let $v_0=v$. For $i\geq0$, let $k_i\in\{0,1,2,3\}$ such that $v_i\in\Sigma_{k_i}$. Define $$v_{i+1}=\sigma^{-1}_{k_i}v_i.$$
If $v$ has slope in $\mathbb Q\bigr[\sqrt 5\bigr]$, then the sequence $\{k_0,k_1,..\}$ is eventually constant equal to $0$, starting from some index $k$. For this value of $k$, $v_k=\big(\begin{smallmatrix}
 l_v\\
 0
 \end{smallmatrix}\big)$ for some $l_v\in\mathbb R$. In other words:
 \begin{theorem}\label{gcd}
\cite{Davis_Lelièvre} Corresponding to any vector $v$ in the first quadrant with slope in $\Q[\sqrt{5}]$ there is a unique word $k_0,k_1,...,k_n$, where $k_i\in\{0,1,2,3\}$, which correspond to an itinerary through the sectors $\Sigma_0,\Sigma_1,\Sigma_2,\Sigma_3$, such that $v=l_v\sigma_{k_0}...\sigma_{k_n} \big(\begin{smallmatrix}
 1\\
 0
  \end{smallmatrix}\big)$, for some $l_v\in\mathbb R_{+}$. The vector in $S$ in this direction is $l_v^{-1}v$.
\end{theorem}
The ring of integers of $\mathbb Q\bigr[\sqrt 5\bigr]$ is $\mathbb Z\bigr[\phi\bigr]$. Since multiplying a pair in $\mathbb Z\bigr[\phi\bigr]^2$ by an element of $\Gamma$ doesn't change the gcd, the iterative process mentioned above is a gcd-algorithm for $\mathbb Z\bigr[\phi\bigr]$. In general, the gcd in $\mathbb Z\bigr[\phi\bigr]$ is defined uniquely only up to multiplication by elements of the unit group. The above algorithm provides a representative for the gcd. In the notations of Theorem \ref{gcd} it is $l_v$.\\  
For $x,y\in \mathbb Z\bigr[\phi\bigr]$, denote by gcd$_{\Gamma^+}(x,y)$ the  representative for the gcd of $|x|,|y|$ obtained by the above algorithm. The following Lemma shows that the $gcd_{\Gamma^+}$ is well defined: \begin{lemma}\label{sym}
    $\big(\begin{smallmatrix}
 x\\
 y
 \end{smallmatrix}\big)\in S \Leftrightarrow \big(\begin{smallmatrix}
 |x|\\
 |y|
 \end{smallmatrix}\big)\in S$.
\end{lemma}
 \begin{proof}
 The orbit S is symmetric with respect to the line $y=x$ (see Lemma 2.17 in \cite{Davis_Lelièvre}). Moreover, $\sigma_0\cdot\sigma_3^{-1}\cdot\sigma_0\cdot\sigma_3^{-1}\cdot\sigma_0=\big(\begin{smallmatrix}
0&-1\\
1&0
 \end{smallmatrix}\big)$ so $\Gamma$ contains a matrix that rotates by $\frac{\pi}{2}$. 
\end{proof}
\lemma{}\label{lemma}
 Let $v_1,v_2  \in S$, let $v_1-v_2=
 \big(\begin{smallmatrix}
 x\\
 y
 \end{smallmatrix}\big)$, and $\lambda=$ gcd$_{\Gamma^+}(x,y)$. Then there are two points in $S$ on a horizontal line in distance $\lambda$.

 \begin{proof}
  $\lambda^{-1}
 \big(\begin{smallmatrix}
 x\\
 y
 \end{smallmatrix}\big) \in S$. So there is $\gamma\in \Gamma$ such that $\gamma \big(\begin{smallmatrix}
 1\\
 0
  \end{smallmatrix}\big) =\lambda^{-1}(v_1-v_2)$, and \[\lambda=\big\|\begin{pmatrix}
  \lambda\\
 0
 \end{pmatrix}\big\|=\big\|\gamma^{-1}v_1-\gamma^{-1}v_2\big\|.\]

 \end{proof}

 \begin{proof}   
 [Proof of Theorem \ref{3points}]For every $k\in\mathbb{N}$, $k>1$ we take 3 vectors in $S$:
 \begin{itemize}
     \item $u_1(k)=\sigma_3^{k+1}\sigma_2
 \big(\begin{smallmatrix}
 1\\
 0
 \end{smallmatrix}\big)= \big(\begin{smallmatrix}
 1 & 0\\
 (k+1)\phi & 1
 \end{smallmatrix}\big)\big(\begin{smallmatrix}
 \phi\\
 \phi
 \end{smallmatrix}\big)=\big(\begin{smallmatrix}
 \phi\\
 k+1+(k+2)\phi
 \end{smallmatrix}\big)$.
 \\
 \item $u_2(k)=\sigma_0^{k-1}\sigma_1\sigma_0^{k-1}\sigma_2
 \big(\begin{smallmatrix}
 1\\
 0
 \end{smallmatrix}\big)= \big(\begin{smallmatrix}
 1 & (k-1)\phi\\
 0 & 1
 \end{smallmatrix}\big) \big(\begin{smallmatrix}
 \phi & \phi\\
 1 & \phi
 \end{smallmatrix}\big) \big(\begin{smallmatrix}
 1 & (k-1)\phi\\
 0 & 1
 \end{smallmatrix}\big)\big(\begin{smallmatrix}
 \phi\\
 \phi
 \end{smallmatrix}\big)$ \[=\big(\begin{smallmatrix}
 k(k+1)+(2k^2+k-1)\phi\\
 k+(k+1)\phi
 \end{smallmatrix}\big).\]
 
 \item $u_3(k)=\sigma_0^{2k+1}\sigma_1\sigma_0^{k-2}\sigma_2
 \big(\begin{smallmatrix}
 1\\
 0
 \end{smallmatrix}\big)= \big(\begin{smallmatrix}
 1 & (2k+1)\phi\\
 0 & 1
 \end{smallmatrix}\big) \big(\begin{smallmatrix}
 \phi & \phi\\
 1 & \phi
 \end{smallmatrix}\big) \big(\begin{smallmatrix}
 1 & (k-2)\phi\\
 0 & 1
 \end{smallmatrix}\big)\big(\begin{smallmatrix}
 \phi\\
 \phi
 \end{smallmatrix}\big)$ \[=\big(\begin{smallmatrix}
 2k(k+1)+(4k^2+2k-3)\phi\\
 k-1+k\phi
 \end{smallmatrix}\big).\]
 \end{itemize}

Computing the difference vectors we get:
\begin{equation}\label{eq:1}
    u_1(k)-u_2(k)=u_2(k)-u_3(k)=
 \begin{pmatrix}
 -(k(k+1)+(2k^2+k-2)\phi)\\
 \phi^2
 \end{pmatrix}
\end{equation} 
  so $u_1(k),u_2(k),u_3(k)$ are on one line. The unit group of the ring $\Z[\phi]$ is generated by $\phi$ and $\overline{\phi}=1-\phi=-\phi^{-1}$. This follows from the fact that $N(\phi)=-1$ together with Dirichlet's unit theorem \cite{units}. Therefore, the coordinates of vector \ref{eq:1} are coprime, and there is an integer $j_k$ such that \begin{equation}\label{eq:2}
  gcd_{\Gamma^+}\bigl(k(k+1)+(2k^2+k-2)\phi,\phi^2\bigl)=\phi^{j_k}.
  \end{equation}
  By Lemma \ref{lemma}, for every  $k>1$ there are 3 elements of $S$ on a horizontal line of length $2\phi^{j_k} $.
 It remains to show that there is a subsequence $\big(j_{k_n}\big)_{n\in\mathbb{N}}$ such that  $\phi^{j_{k_n}}\xrightarrow[n \to \infty]{}0$. 
 \\
 
 Denote by $d_k$ the difference vector reflected into the first quadrant:
\begin{equation}\label{dk}d_k=\begin{pmatrix}
 k(k+1)+(2k^2+k-2)\phi\\
 \phi^2
 \end{pmatrix}=\begin{pmatrix}
 k^2\phi^3+k\phi^2-2\phi\\
 \phi^2
 \end{pmatrix}
\end{equation} where the second equality follows from $\phi$'s property: $\phi^2=1+\phi, \hspace{1mm}\phi^3=1+2\phi$. Recall that in the $gcd_{\Gamma^+}$-algorithm described in Theorem \ref{gcd}, if $v_i\in\Sigma_{k_i}$, $$v_{i+1}=\sigma^{-1}_{k_i}v_i.$$ Moreover, it  follows from multiplying a vector by the matrices $\sigma_0, \sigma_1,\sigma_2,\sigma_3$, that for every $v$ in the first quadrant and every $j\in\{0,1,2,3\}$, $$\Pi_x(v)\leq\Pi_x(\sigma_jv)$$ where $\Pi_x$ is the projection on the first axis. Hence, for every $i$, $$\Pi_x(v_{i+1})\leq \Pi_x(v_i).$$ Clearly for every $k$, $d_k\in \Sigma_0$. Thus, the first step of the $gcd_{\Gamma^+}$ algorithm on $d_k$, is multiply by $\sigma_0^{-1}$. Let $t\geq1$ be the integer such that $\sigma_0^{-(t-1)}d_k\in\Sigma_0$ and $\sigma_0^{-t}d_k\notin\Sigma_0$ ($t$ exist by Theorem \ref{gcd}). In order to find $t$, we compute $t$ such that $\sigma_0^{-t}d_k$ is on the line $y=\phi^{-1}x$, and then take the ceiling of this value. We get: \[t=\Bigl\lfloor \frac{k^2\phi^3+k\phi^2-2\phi}{\phi^3}\Bigr\rfloor=\bigl\lfloor k^2+k(\phi -1)-2(2-\phi)\bigr\rfloor= k^2-k-4+\bigl\lfloor(k+2)\phi\bigr\rfloor.\] 
 Since the $gcd_{\Gamma^+}$ is attained by $\Pi_x(v_i)$ when $\Pi_y(v_i)=0$, we conclude:
 $$\phi^{j_k}=gcd_{\Gamma^+}(d_k)\leq
\Pi_x(\sigma_0^{-t}d_k)=
k^2\phi^3+k\phi^2-2\phi-(k^2-k-4+\bigl\lfloor(k+2)\phi\bigr\rfloor)\phi^3$$
$$=\phi^3(k+4-\bigl\lfloor(k+2)\phi\bigr\rfloor)+k\phi^2-2\phi$$
$$=\phi^3(k+4+k\phi^{-1}-2\phi^{-2}-\bigl\lfloor(k+2)\phi\bigr\rfloor)=\phi^3((k+2)\phi-\bigl\lfloor(k+2)\phi\bigr\rfloor)$$
$$=\phi^3\bigl\{(k+2)\phi\bigr\}.$$
\\
If $k+2=F_n$, where $n$ is odd, we have $$\phi^{j_k}\leq \phi^3\bigl\{(k+2)\phi\bigr\}=\phi^3\{F_n\phi\}=\phi^3(F_n\phi-\lfloor F_n\phi\rfloor)=\phi^3\phi^{-n}.$$
Hence, there is a subset  $(k_n)_{n\in\mathbb{N}}\subseteq\mathbb {N}$ such that $$\phi^{j_{k_n}}\xrightarrow[n \to \infty]{}0.$$
 $$$$
 We now estimate the norm of the points. For $k\in \N$, using the notations of Lemma \ref{lemma}, there is $\gamma\in \Gamma$ such that the points can be presented as $\gamma^{-1}u_i(k)$, and 
$$\gamma(\begin{smallmatrix}
    1\\0 \end{smallmatrix})=\phi^{-j_k}(u_1(k)-u_2(k))=\phi^{-j_k}(\begin{smallmatrix}
 1 & 0\\
 0 & -1
 \end{smallmatrix})\cdot d_k=\phi^{-j_k}
\bigl(\begin{smallmatrix}
 k^2\phi^3+k\phi^2-2\phi\\
 -\phi^2
 \end{smallmatrix}\bigr).$$ As $det(\gamma)=1$, we invert $\gamma$ and find that the second row of $\gamma^{-1}$ is $$\phi^{-j_k}\begin{pmatrix}
 \phi^2\\
 k^2\phi^3+k\phi^2-2\phi
 \end{pmatrix}^T.$$
 Notice that $\gamma^{-1}$ takes the difference vector in equation (\ref{eq:1}) to a horizontal line. Hence, $\gamma^{-1} u_1(k),\gamma^{-1} u_2(k),\gamma^{-1} u_3(k)$ are within a segment of length $2\phi^{j_k}$ which has the height: $$ \phi^{-j_k}\begin{pmatrix}
 \phi^2\\
 k^2\phi^3+k\phi^2-2\phi
 \end{pmatrix}\cdot u_1(k)=\phi^{-j_k}\begin{pmatrix}
 \phi^2\\
 k^2\phi^3+k\phi^2-2\phi
 \end{pmatrix}^T \begin{pmatrix}
 \phi\\
 (k+1)\phi^2+\phi
 \end{pmatrix}$$$$=\phi^{-j_k}\bigl(\phi^3+k^2(k+1)\phi^5+k^2\phi^4+k(k+1)\phi^4+k\phi^3-2(k+1)\phi^3-2\phi^2\bigr)$$$$=\phi^{-j_k}\bigl(k^3\phi^5+k^2(\phi^4+\phi^6)-(\phi^2+\phi^4)\bigr).$$ We can assume that these points are not in $\Sigma_0$ (otherwise act with $\sigma_0^{-1})$. So, for $i=1,2,3$ $$\big\|\gamma^{-1}u_i(k)\big\|\leq\big\| \gamma^{-1}u_3(k)\big\|\leq c \cdot \phi^{-j_k}\bigl(k^3\phi^5+k^2(\phi^4+\phi^6)-(\phi^2+\phi^4)\bigr)
   =O(\phi^{-j_k}k^3).$$ 
\\
Again, we take $k+2=F_n$, where $n$ is odd. We have $$\phi^{j_k}\leq \phi^3\bigl\{(k+2)\phi\bigr\}=\phi^3\{F_n\phi\}=\phi^3(F_n\phi-\lfloor F_n\phi\rfloor)=\phi^3\phi^{-n}.$$
Hence, $k=\frac{1}{\sqrt5}(\phi^n-(-\phi)^{-n})-2=O(\phi^n)=O(\phi^{-j_k}).$
Finally, $$\big\|\gamma^{-1}u_i(k)\big\|=O(\phi^{-j_k}k^3)=O((\phi^{j_k})^{-4}).$$
Thus, $\epsilon_n=2\cdot\phi^{j_{(F_n-2)}}$ is the desired sequence.
\end{proof}
\normalfont
We give an alternate proof of the first part of Theorem \ref{3points}. With further work, this approach may lead to improved quantitative results or offer insights into the orbits of other non-uniform lattices.
\begin{proposition}\label{remark}
Let $\big(j_k\big)$ as in equation (\ref{eq:2}) in the preceding proof. For every $m\in \Z$, there is a finite number of indices $k$ such that $j_k=m$.
\end{proposition}

\begin{proof}
Every point of $S$ in the first quadrant can be written as $\sigma_0^{t}z$ where $t\in\mathbb{N}$, and $z\in S':=S\cap\bigl(\Sigma_1\cup\Sigma_2\cup\Sigma_3\bigl)$. In particular for every $k>1$ there is $z\in S'$ and $t\in \mathbb{N}$ such that $$\phi^{-j_k}d_k=\sigma_0^{t}z.$$ 
Assume to the contrary that there exists an $m$, and infinitely many $k's$ such that $j_k=m$. The shortest vectors corresponding to words of length $n$ are $(\begin{smallmatrix}
    1\\n\phi \end{smallmatrix}),(\begin{smallmatrix}
    n\phi\\1 \end{smallmatrix})$ (see Proposition 2.13 in \cite{Davis_Lelièvre}). Since there are finitely many words of length less than $n$, it follows that every bounded set in the first quadrant contains a finite number of points from $S$. In particular the set $S'\cap\{y=\phi^{2-m}\}$  is finite because it is bounded. But from the assumption there are infinitely many points of the form $\phi^{-m}d_k$ on the positive ray along the line $y=\phi^{2-m}$.
Hence, there is $z\in S'$ and infinite subset $\big(k_n\big)\subset \N$ such that  $$\phi^{-m}d_{k_n}=\sigma_0^{t_{k_n}} z.$$
The height of $d_{k_n}$ is $\phi^2$ (\ref{dk}). Therefore, the distance between any two points of the form given by the right-hand side of the equation is an integer multiple of $\phi^{3-m}$.
In other words, for every $n_1,n_2\in \mathbb{N}$, there is $t\in\mathbb{Z} $ such that 
\[\phi^{-m}\big(d_{k_{n_1}}-d_{k_{n_2}}\big)=\begin{pmatrix}
    0 & (t_{k_{n_1}}-t_{k_{n_2}})\phi \\
    0 & 0
\end{pmatrix}\cdot z = t\begin{pmatrix}
 \phi^{3-m}\\
 0
 \end{pmatrix}.\] i.e.,
 \[\phi^{-m}\begin{pmatrix}
 (k_{n_1}^2-k_{n_2}^2)\phi^3+(k_{n_1}-k_{n_2})\phi^2\\
 0
 \end{pmatrix}=t\begin{pmatrix}
 \phi^{3-m}\\
 0
 \end{pmatrix}.\]
 By multiplying both parts by $\phi^{m-3}$ we get: \[\bigr(k_{n_1}^2-k_{n_2}^2\bigr)+\bigr(k_{n_1}-k_{n_2}\bigr)\phi^{-1}=t\]
 \[\Rightarrow \bigr(k_{n_1}-k_{n_2}\bigr)\phi^{-1}=t- \bigr(k_{n_1}^2-k_{n_2}^2\bigr) \in \Z.\]
 \
 Therefore $k_{n_1}=k_{n_2}$,
 a contradiction to $\big|\big(d_{k_n}\big)\big|=\infty$.\\
 \end{proof}

\begin{proof}[Alternate proof of first part of Theorem \ref{3points}]
Since $(\begin{smallmatrix}
    1\\0 \end{smallmatrix})$ is the shortest vector of $S$ in the first quadrant, there are only finitely many indices $k$ such that $j_k>0$. It follows from Proposition \ref{remark} that $\lim_{k\to\infty} (j_k)=-\infty$.
\end{proof}
\section*{Questions}
\begin{enumerate}
\item Is S relatively dense? i.e., is there $R>0$ such that any ball of radius $R$ contains at least one point of $S$? For rational $\rho$, it is known that for $\Gamma_\rho
\mathbf{e}_1$ as in Proposition \ref{rho}, the answer is no \cite{Wu}.
\item What is the decreasing rate of the sequence $\big(j_k\big)$ as in equation (\ref{eq:2})?
\item Does Conjecture \ref{conjecture} hold for orbits of other non-uniform lattices in $SL_2(\R)$, particularly those arising from other non-arithmetic Veech surfaces?
\end{enumerate}

 \end{document}